\documentclass{amsart}  %
\catcode`\@=11

 \def\AMSTeXfeatures{\Plainheads 
   \let\current@vert=\AMS@vert}

 \def\Plainheads{\sh@ftdiam=0.05em
   \getlabeldims
   \let\vshaftfill=\plnvsolidfill
   \let\hshaftfill=\plnhsolidfill
   \let\th@rhead=\plnrhead
   \let\th@lhead=\plnlhead
   \let\th@dnhead=\plndnhead
   \let\th@uphead=\plnuphead}
 
 \def\glet{\global\let}

 \def\LaTeXfeatures{\catcode`\@=11
   \ifx\@clnwd\undefined \nol@g
      \input ltxcode.tex \dol@g \fi
   \ltxheads \let\current@vert=\new@vert
   \providelto \catcode`\@=\active}

 \def\nol@g{\def\wlog{\edef\garbage}}
 \def\dol@g{\let\wlog=\wl@g} \let\wl@g=\wlog
 \nol@g % Silences allocations with \newbox, etc.

 \newbox\ltobox
 \def\providelto{{\setbox\z@=
   \hbox{$\to$}\minharrlen=\wd\z@
   \global\setbox\ltobox=\hbox{$\activeat>>>$}}
   \def\lto{\mathrel{\copy\ltobox}}}

 \def\ltxheads{\sh@ftdiam=\@wholewidth
   \getlabeldims
   \let\vshaftfill= \ltxvsolidfill
   \let\hshaftfill=\ltxhsolidfill
   \let\th@rhead=\ltxrhead
   \let\th@lhead=\ltxlhead
   \let\th@dnhead=\ltxdnhead
   \let\th@uphead=\ltxuphead}
 {\catcode`\@=\active
   \gdef@#1{\csname #1\string@at\endcsname}
   \glet\activeat=@}
 \def\def@#1{\expandafter\def\csname #1@at\endcsname}

 \def@>#1>#2>{\@rrow R{#1}{#2}}
 \def@<#1<#2<{\@rrow L{#1}{#2}}
 \def@ V#1V#2V{\@rrow V{#1}{#2}}
 \def@ A#1A#2A{\@rrow A{#1}{#2}}
 \def@/#1/#2/#3/{\@rrow{#1}{#2}{#3}}
 \def@.{\ifodd\row\ifmmode\noharrow
     \else\leavevmode.\spacefactor3000 \fi
   \else\novarrow\fi}
 \def@={\ifodd\row\harrow\hequalfill{}{}%
   \else\varrow\vequalfill{}{}\fi}
 \def@:#1{\ifx=#1\harrow\deffill{}{}%
   \else\leavevmode\null:#1\fi}
 \def@|{\current@vert}
  \def\AMS@vert{\varrow\vequalfill{}{}}
  \def\new@vert#1|#2|{\ifodd\row
   \let\nextarrow\vertexvarrow
   \else\let\nextarrow\varrow\fi
   \nextarrow\vshaftfill{#1}{#2}}
 \def@-{\ifmmode\let\next\hl@ne
   \else\let\next\AMSatdash \fi \next}
  \def\hl@ne#1-#2-{\harrow\hshaftfill{#1}{#2}}
  \def\AMSatdash{\let\next\relax\leavevmode
    \def\next@{\ifx\next-%
      \def\next-{\futurelet\next\nextii@}%
     \else\def\next{\hbox{-}}\fi\next}%
    \def\nextii@{\ifx\next-\def\next-{\hbox{---}}%
      \else\def\next{\hbox{--}}\fi\next}%
    \futurelet\next\next@}
 \def@(#1){\tweenarrows{#1}}
 \def@[#1]{\setsp@n#1\relax\activeat}
 \def\fiberbox{\hbox{$\vcenter{\hr@le\hbox{\vr@le
   \kern1ex\vbox{\kern1.2ex}\vr@le}\hr@le}$}}
  \def\hr@le{\hrule height \sh@ftdiam}
  \def\vr@le{\vrule width \sh@ftdiam}
 \def@+#1+#2+#3+{\ifodd\row \harrow{#1}{#2}{#3}%
   \else \varrow{#1}{#2}{#3}\fi}

 \def\Dnarrfill{\vequalfill\Dnhe@d}
 \def\Uparrfill{\Uphe@d\vequalfill}
 \def\ontofill{\rtarrfill\kern-0.3em %2\he@dwd
   \th@rhead\kern 0.3em} %new def

 \def\rtarrfill{\hshaftfill\th@rhead}
 \def\ltarrfill{\th@lhead\hshaftfill}
 \def\dnarrfill{\vshaftfill\th@dnhead}
 \def\uparrfill{\th@uphead\vshaftfill}
 \def\hequalfill{\plnhfill=}
 \def\deffill{:\plnhfill=}
 \def\plnvextfill#1{\setbox\z@
   \hbox{\the\textfont3 #1}%
   \dimen@=\dp\z@\advance\dimen@\ht\z@
   \copy\z@ \kern-\dimen@ %-\dp\z@
   \cleaders\copy\z@ \vfill
   \kern-\dimen@ %-\dp\z@
   \box\z@}
 \def\plnhfill#1{$\m@th\mkern-1.5mu\mathord#1\mkern-6mu
    \cleaders\hbox{$\mkern-2mu\mathord#1\mkern-2mu$}\hfill
    \mkern-6mu\mathord#1\mkern-1.5mu$}
 \def\vequalfill{\plnvextfill{\char'167}}
 \def\plnvsolidfill{\plnvextfill{\char'077}}
 \def\plnhsolidfill{\plnhfill-}
 \def\ltxhsolidfill{\leaders\hrule height\topofshaft depth\botofshaft
   \hfill}
 \def\ltxvsolidfill{\leaders\vrule width\sh@ftdiam\vfill}
 \def\hdashfill{\hd@sh\wd@sh
   \xleaders \hbox{\wd@sh\hd@sh\wd@sh}\hfill
   \wd@sh\hd@sh}
 \def\vdashfill{\vd@sh\wd@sh
   \xleaders \vbox{\wd@sh\vd@sh\wd@sh}\vfill
   \wd@sh\vd@sh}
 \def\dashed{\ifinmeasureCD\else
    \ifodd\row\option{\let\hshaftfill=\hdashfill}%
   \else\option{\let\vshaftfill=\vdashfill}\fi\fi}
 \newdimen\CDstrutht  \newdimen\CDstrutdp
 \newdimen\CDstrutlen \CDstrutlen=\CDstrutht
   \advance\CDstrutlen by \CDstrutdp
 \def\CDstrut{\vrule
   height \ifnum\row=1 \z@\else\CDstrutht \fi
   depth \ifnum\row=\numrows \z@ \else\CDstrutdp \fi
   width\z@}
 \newdimen\CDarrsurr \CDarrsurr=0.375em
 \newdimen\CDdashlen
 \newdimen\CDvarrlen \CDvarrlen=1.5\baselineskip
 \newdimen\minharrlen %Used outside CD's
  \setbox\z@\hbox{$\longrightarrow$} \minharrlen=\wd\z@
 \newdimen\minCDharrlen \minCDharrlen=2.5em %825 2pc %2.5pc
\newdimen \minc@lwd
\def\findminc@lwd{\minc@lwd=2\CDarrsurr
  \advance\minc@lwd\minCDharrlen}
 \newdimen\sh@ftdiam

 \newdimen\labelsurr \labelsurr=1.25 em

\newcount\sp@ncnt \sp@ncnt=\@ne
\newcount\sp@ncnt@ \sp@ncnt@=\@ne
\newdimen\@rrwd \newdimen\@rrdp

 \def\adjustbot#1{\option{\advance\@rrdp#1\relax}}
\def\pushvertex#1{\global\p@shlen#1\relax
   \global\let\maybepush=\dopush}

 \newdimen\p@shlen \p@shlen=\z@

 \let\maybepush=\relax
 \def\dopush{\ifinmeasureCD %omitted by accident
   \advance\locdimen by -\p@shlen %AL
   \else\advance \@rrwd by -\p@shlen \fi %AL
   \global\let\maybepush=\relax \global\p@shlen=\z@\relax}
 \def\span@ne{\global\sp@ncnt=\@ne\relax}
 \def\setsp@n#1#2{\global\sp@ncnt=#1\relax
   \ifx\relax#2\relax\else\global\sp@ncnt@=#2\relax\fi}

 \def\plnrhead{\llap{$\rightarrow\mkern-1.5mu$}}
 \def\plnlhead{\rlap{$\mkern-1.5mu\leftarrow$}}

 \def\clap#1{\hbox to \z@{\hss #1\hss}}

 \def\plndnhead{\hbox{\the\textfont3 \char'171}}
 \def\plnuphead{\hbox{\the\textfont3 \char'170}}
 \def\Dnhe@d{\hbox{\the\textfont3 \char'177}}
 \def\Uphe@d{\hbox{\the\textfont3 \char'176}}

 \def\ltxrhead{\raise\@xisheight
   \llap{\smash{\@linefnt\@getrarrow(1,0)}}}
 \def\ltxlhead{\raise\@xisheight
   \rlap{\@linefnt\@getlarrow(-1,0)}}
 \def\ltxuphead{\setbox\z@=\rlap{%
   \kern\@halfwidth\@linefnt\char'66}%
   \copy\z@\kern-\ht\z@}
 \def\ltxdnhead{\setbox\z@=\rlap{%
   \kern\@halfwidth\@linefnt\char'77}%
   \ht\z@=\z@\box\z@}

 \def\wd@sh{\kern0.5\CDdashlen}
 \def\hd@sh{\vrule height\topofshaft depth\botofshaft
    width\CDdashlen}
 \def\vd@sh{\hrule height\CDdashlen
   depth\z@ width\sh@ftdiam}
\def\xylist{14{3434}13{2414}12{1723}%
  23{1413}34{1153}11{0867}43{0707}%
  32{0580}21{0414}31{0291}41{0}}
\newcount\tgtcnt@
\def\find@xyargs{\dimen@=\@rrdp
  \advance\dimen@ by \CDstrutlen
  \tgtcnt@=\dimen@ \dimen@=\@rrwd %\relax
  \divide\dimen@ by \@m %\relax
  \divide \tgtcnt@ by \dimen@ %\relax
  \expandafter\testxy\xylist\relax
  \unitlength=\@xarg\@rrdp
  \divide\unitlength by\@yarg\relax}
\def\testxy#1#2#3{\ifnum\tgtcnt@>#3
    \@xarg=#1\relax \@yarg=#2\relax
    \let\next=\ignorerest
  \else\let\next\testxy\fi\next}
\def\ignorerest#1\relax{\relax}

\let\scalefactor=\@ne
\def\SWarrow{\find@xyargs\vector
  (-\@xarg,-\@yarg)\scalefactor\hskip-\wd\@linechar}
\def\NWarrow{\find@xyargs\vector
  (-\@xarg,\@yarg)\scalefactor\hskip-\wd\@linechar}
\def\NEarrow{\find@xyargs\vector
  (\@xarg,\@yarg)\scalefactor}
\def\SEarrow{\find@xyargs\vector
  (\@xarg,-\@yarg)\scalefactor}
\def\rightupline{\find@xyargs\@linelen=\scalefactor
     \unitlength\@sline}
\def\rightdownline{\find@xyargs\@yarg=-\@yarg\relax
     \@linelen=\scalefactor\unitlength\@sline}

\def\Sim{\ifodd\row\setbox\z@=\hbox{$\sim$}\dimen@=\ht\z@
 \advance\dimen@ by -\@xisheight
  \vbox{\box\z@\kern-\@xisheight\kern\dimen@}%
  \else\hbox{$\wr$}\fi}

\def\harrow#1#2#3{\inmeasureCDtrue\findminarrwd
  {#2}{#3}{\sp@ncnt\minharrlen}\inmeasureCDfalse\span@ne
  \mathrel{\hbox{\options\hplace{#1}\ulabel{#2}\dlabel{#3}}}}

\def\noharrow{\harrow\hfill{}{}}
\def\vertexvarrow#1#2#3{\findarrdp \@rrwd=\z@ \setsp@n\@ne\@ne
  \vbox to \z@{\kern-1.2\CDstrutht
  \rlap{\options\vplace{#1}\llabel{#2}\rlabel{#3}}\vss}}

\newif\ifinmeasureCD
\def\measurelabel#1{\setbox\z@
  \hbox{$\scriptstyle#1\kern\labelsurr$}%
  \ifdim\wd\z@>\@rrwd \@rrwd=\wd\z@\fi}
\def\findminarrwd#1#2#3{\@rrwd=#3\relax
   \measurelabel{#1}\measurelabel{#2}}
\def\findCDarrwd#1#2{\@rrwd=\minCDharrlen
   \measurelabel{#1}\measurelabel{#2}%
  }

\newcount\row \row=\@ne \newcount\col \col=\@ne %96 initialized
 \newcount\numrows 
\numrows=\@ne
 \newcount\numcols
\newcount\arrspan \newdimen\vrtxhalfwd  \newbox\tempbox

\def\DANABUG{\advance\col by \@ne
 \@rrwd=\minCDharrlen
  \advance\@rrwd by \vrtxhalfwd
  \advance\@rrwd by \CDarrsurr
  \ifnum\col>\numcols \numcols=\col
     \newlocdimen{col\the\col}\locdimen=\@rrwd %AL
  \else \ifdim\@rrwd>\c@l \c@l=\@rrwd\fi\fi}

\def\drop#1\\{%\noharrow %caused by DANABUG
  \findvrtxhalfsum\DANABUG\advance\row by 2 \measureinit}

\def\measureinit{\col=\@ne \vrtxhalfwd=-\CDarrsurr\arrspan=\@ne\@rrwd=\z@
   \setbox\tempbox=\hbox\bgroup$}
\def\measure{%CR \bgroup
  \let\harrow\measureCDarrow
  \let\CDCR=\measureCR %CR
   \findminc@lwd 
  \inmeasureCDtrue
  \row=\@ne \numcols=\z@ \measureinit}

\def\endmeasure{\findvrtxhalfsum\DANABUG
  \numrows=\row %CR \egroup
  \inmeasureCDfalse}

\def\newlocdimen#1{\advance\dimenc@unt by \@ne
  \ifnum\dimenc@unt<\insc@unt
     \else\errmessage{No room for the CD}\fi
  \dimendef\locdimen=\dimenc@unt
  \expandafter\dimendef\csname#1\endcsname=\dimenc@unt}

 \def\r@wc@l{\csname row\the\row col\the\col\endcsname}
 \def\c@l{\csname col\the\col\endcsname}

 \def\findvrtxhalfsum{$\egroup
  \newlocdimen{row\the\row col\the\col}%%AL
  \locdimen=\vrtxhalfwd %AL
  \vrtxhalfwd=0.5\wd\tempbox %\maybes@ve %8231
  \advance\vrtxhalfwd by \CDarrsurr
  \advance\locdimen by \vrtxhalfwd %AL
  \advance\@rrwd by \locdimen %AL
  \maybepush
  \divide\@rrwd by \arrspan\relax
  \ifdim\@rrwd<\minc@lwd
    \ifnum\col>\@ne \@rrwd=\minc@lwd\fi \fi
  \loop %94 \edef\c@l{\csname col\the\col\endcsname}
    \ifnum\col>\numcols \numcols=\col
       \newlocdimen{col\the\col}% %AL
       \locdimen=\@rrwd %AL
    \else \ifdim\@rrwd>\c@l \c@l=\@rrwd\fi \fi
   \ifnum\arrspan>\@ne
      \advance\arrspan by -1 \advance\col by \@ne
  \repeat }

 \def\measureCDarrow#1#2#3{\findvrtxhalfsum
   \arrspan=\sp@ncnt\relax\global\sp@ncnt=1\relax
   \advance\col by \@ne
   \findCDarrwd{#2}{#3}%
   \setbox\tempbox=\hbox\bgroup$}

 \newcount\dr@tn \dr@tn=\z@
 \def\locate#1:#2{\ifinmeasureCD\else
   \count@=-#1
   \multiply\count@ by 2
   \advance\count@ by #2
   \dimen@=\count@\@rrwd
   \ifnum\dr@tn=\@ne\relax \else\dimen@=-\dimen@ \fi
   \dimen@i=\@rrdp
   \ifnum\dr@tn>\z@\advance\dimen@i by \CDstrutlen \fi
   \dimen@i=\count@\dimen@i
   \count@=#2 \multiply\count@ by 2
   \divide\dimen@ by \count@
   \divide\dimen@i by \count@
   \lift\dimen@i\nudge\dimen@\fi}
\def\betweenCDrows{\advance\row by \@ne \col=\@ne
\options}

\def\hbegin{\hbox\bgroup\kern\c@l \kern-\r@wc@l$}
\def\hend{$\glet\maybepush\relax \CDstrut\egroup}
\def\vbegin{\setbox\tempbox=\hbox\bgroup$}
\def\vend{$\egroup\ht\tempbox=\z@\dp\tempbox\CDvarrlen
  \box\tempbox}
\def\setCD{\let\harrow=\setCDarrow
  \let\CDCR=\setCR %CR
  \row=\@ne \col=\@ne \hbegin}
\let\endsetCD=\hend %AT (Assume CD ends with hmaterial)

\def\findarrwd{\@rrwd=\z@ \count@=\col \advance\count@ by\sp@ncnt
  \loop\ifnum\count@>\col \advance\count@ by -1
      \advance\@rrwd by\csname col\the\count@\endcsname\repeat}
\def\setCDarrow#1#2#3{\kern\CDarrsurr\advance\col by \@ne
  \findarrwd \advance\@rrwd by -\r@wc@l  
  \@rrdp=\z@ %&0211 (It might be used by \locate).
  \maybepush
  \advance\col by -\@ne \advance\col by \sp@ncnt \span@ne
  \hbox to \@rrwd{\options
   \@rrwd=\scalefactor\@rrwd\hss
   \hplace{#1}\ulabel{#2}\dlabel{#3}\hss}%
   \kern\CDarrsurr}

\newdimen\labspacei %96 use subscript min of TeXbook 13a, p.444
\newdimen\labspaceii %96 Note many letters stick down below their boxes.

\newdimen\@xisheight
  \@xisheight=\the\fontdimen22\textfont2
\newdimen\labelskip
  \labelskip=\the\fontdimen10\textfont3 %2pt
\newdimen\topofshaft
\newdimen\botofshaft
\newdimen\botofulabel
\newdimen\topofdlabel
\def\getlabeldims{
  \topofshaft=0.5\sh@ftdiam
  \botofshaft=\topofshaft
  \advance\topofshaft by \@xisheight  
  \advance\botofshaft by -\@xisheight  
  \botofulabel=\topofshaft
  \advance\botofulabel by \labelskip
  \topofdlabel=\botofshaft
  \advance\topofdlabel by \labelskip}

\def\ulabel{\ifnum\row=\@ne\let\next\ulabeli
   \else\let\next\ulabellap\fi\next}
\def\ulabeli#1{\vbox{
  \clap{\kern-\@rrwd$\scriptstyle#1$}%
  \kern\botofulabel}\maybeoffset}
\def\ulabellap#1{\vbox to \z@{\vss
  \clap{\kern-\@rrwd$\scriptstyle#1$}%
  \kern\botofulabel}\maybeoffset}
\def\dlabel{\ifnum\row=\numrows\let\next\dlabeli
   \else\let\next\dlabellap\fi\next}
\def\dlabeli#1{\vtop{\kern\topofdlabel
  \clap{\kern-\@rrwd$\scriptstyle#1$}%
  }\maybeoffset}
\def\dlabellap#1{\vbox to \z@{\kern\topofdlabel
  \clap{\kern-\@rrwd$\scriptstyle#1$}%
  \vss}\maybeoffset}
\def\rlabel#1{\vbox to \z@{\vss
  \rlap{\kern\labelskip$\scriptstyle#1$}%
  \vss\kern-\@rrdp}\maybeoffset}
\def\llabel#1{\vbox to \z@{\vss
  \llap{$\scriptstyle#1$\kern\labelskip}%
  \vss\kern-\@rrdp}\maybeoffset}
\def\swlabel#1{\vtop{\kern0.5\@rrdp
  \llap{$\scriptstyle#1$\kern\labelskip\kern-0.5\@rrwd}
  }\maybeoffset}
\def\nwlabel#1{\vbox{
  \llap{$\scriptstyle#1$\kern\labelskip\kern-0.5\@rrwd}%
  \kern-0.5\@rrdp}\maybeoffset}
\def\selabel#1{\vtop{\kern0.5\@rrdp
  \rlap{\kern0.5\@rrwd\kern\labelskip$\scriptstyle#1$}%
  }\maybeoffset}
\def\nelabel#1{\vbox{
  \rlap{\kern0.5\@rrwd\kern\labelskip$\scriptstyle#1$}%
  \kern-0.5\@rrdp}\maybeoffset}
\def\cplace#1{\vbox to \z@{\vss
  \clap{$#1$\kern-\@rrwd}%
  \kern-\@rrdp\vss}\maybeoffset}
\def\hplace#1{\hbox to \@rrwd{#1}\maybeoffset}
\def\vplace#1{\clap{\vbox to \z@{#1\kern-\@rrdp}}\maybeoffset}
\newdimen\nudgeamount \nudgeamount=\z@
\newdimen\liftamount \liftamount=\z@
\let\maybeoffset\relax
\newbox\offsetbox \newdimen\lastheight
\def\dooffset{%assumes that \lastbox is a <box> set in horiz. mode
  \setbox\offsetbox=\lastbox \lastheight=\ht\offsetbox 
  \setbox\offsetbox=\vbox{\kern-\liftamount\box\offsetbox}%
  \ht\offsetbox=\lastheight
  \kern\nudgeamount\box\offsetbox\kern-\nudgeamount
  \global\nudgeamount=\z@ \global\liftamount=\z@
  \glet\maybeoffset=\relax}
\def\nudge#1{\ifinmeasureCD\else
  \global\advance\nudgeamount#1\relax
  \global\let\maybeoffset\dooffset\fi}
\def\lift#1{\ifinmeasureCD\else
  \global\advance\liftamount#1\relax
  \global\let\maybeoffset\dooffset\fi}

\def\findarrdp{\@rrdp=\CDvarrlen
  \ifnum\sp@ncnt@>1
    \advance\@rrdp by \CDstrutlen
    \multiply\@rrdp by \sp@ncnt@
    \advance\@rrdp by -\CDstrutlen \fi
 }

\def\varrow#1#2#3{\ifnum\sp@ncnt>\@ne 
     \sp@ncnt@=\sp@ncnt\relax\fi
  \findarrdp \@rrwd=\z@ %&0211 It might be used by \locate
  \kern\c@l
   \hbox to \z@{\options
   \@rrdp=\scalefactor\@rrdp
    \hss\vplace{#1}\llabel{#2}\rlabel{#3}\hss}%
  \global\advance\col by \@ne \setsp@n\@ne\@ne
  }

\def\novarrow{\varrow\vfill{}{}}

\def\tweenarrows#1{\findarrwd \findarrdp \setsp@n\@ne\@ne
  \rlap{\options\cplace{#1}}}

\def\usarrow #1#2#3{\dr@tn=\@ne
  \findarrwd \findarrdp \setsp@n\@ne\@ne 
  \rlap{\options\cplace{#1}\nwlabel{#2}\selabel{#3}}%
  \dr@tn=\z@}
\def\dsarrow #1#2#3{\dr@tn=\tw@
  \findarrwd \findarrdp \setsp@n\@ne\@ne 
  \rlap{\options\cplace{#1}\swlabel{#2}\nelabel{#3}}%
  \dr@tn=\z@}
 \def\@rrow#1{\csname #1@rrow\endcsname}
 \def\R@rrow{\harrow \rtarrfill}
 \def\L@rrow{\harrow \ltarrfill}
 \def\V@rrow{\varrow \dnarrfill}
 \def\A@rrow{\varrow \uparrfill}
 \def\SE@rrow{\dsarrow \SEarrow}
 \def\NW@rrow{\dsarrow \NWarrow}
 \def\SW@rrow{\usarrow \SWarrow}
 \def\NE@rrow{\usarrow \NEarrow}
 \def\DS@rrow{\dsarrow \dnslope}
 \def\US@rrow{\usarrow \upslope}
 \def\upslope{\find@xyargs
       \@linelen=\unitlength\@sline}
 \def\dnslope{\find@xyargs\@yarg=-\@yarg\relax
       \@linelen=\unitlength\@sline}

\newtoks\optionlist 
\optionlist={}
\let\options\relax
\def\dooptions{\the\optionlist\global\optionlist={}%
  \glet\options=\relax}
\def\option#1{\ifinmeasureCD\else
  \glet\options=\dooptions
  \global\optionlist=\expandafter{\the\optionlist\relax#1}\fi}
\def\wider#1{\ifinmeasureCD\else
   \option{\advance\@rrwd by #1}\fi}
\def\deeper#1{\ifinmeasureCD\else
   \option{\advance\@rrdp by #1}\fi}
{\def\\{\global\let\sptoken= }\\ }%now \sptoken is a spacetoken

\def\CR{\futurelet\nexttok\testCR}
\def\testCR{\ifx\nexttok\sptoken
   \let\next\eatspaceCR\else\let\next\CDCR\fi\next}
\def\eatspaceCR#1 {\CR}
\def\measureCR{\ifx\nexttok\endmeasure\let\nextCR\relax
    \else\let\nextCR\drop\fi\nextCR}
\def\setCR{\ifodd\row
  \ifx\nexttok\endsetCD\else\hend\betweenCDrows\vbegin\fi
  \else\vend\betweenCDrows\hbegin\fi}

\countdef\dimenc@unt=11
\def\CD#1\endCD{%CRAL
   \begingroup\let\\=\CR
  \m@th\offinterlineskip
   \measure#1\endmeasure\null\,\vcenter{\setCD#1\endsetCD}\,
   \endgroup
    }

\ifx\@clnwd\undefined \nol@g\else\catcode`\ =14\relax\fi
 \font\@linefnt=line10 
 \newcount\@tempcnta
 \newcount\@tempcntb
 \newdimen\@tempdima
 \newdimen\@tempdimb
 \newdimen\@wholewidth
 \newdimen\@halfwidth
   \@wholewidth\fontdimen8\@linefnt \@halfwidth .5\@wholewidth
 \newdimen\unitlength
 \newcount\@xarg
 \newcount\@yarg
 \newcount\@yyarg
 \newbox\@linechar
 \newdimen\@linelen
 \newdimen\@clnwd
 \newdimen\@clnht
 \newif\if@negarg
 
 \def\@whilenoop#1{}

 \def\@whiledim#1\do #2{\ifdim #1\relax#2\@iwhiledim{#1\relax#2}\fi}

 \def\@iwhiledim#1{\ifdim #1\let\@nextwhile=\@iwhiledim 
         \else\let\@nextwhile=\@whilenoop\fi\@nextwhile{#1}}

 \def\@sline{\ifnum\@xarg< 0 \@negargtrue \@xarg -\@xarg \@yyarg -\@yarg
   \else \@negargfalse \@yyarg \@yarg \fi
 \ifnum \@yyarg >0 \@tempcnta\@yyarg \else \@tempcnta -\@yyarg \fi
 \ifnum\@tempcnta>6 \@badlinearg\@tempcnta0 \fi
 \ifnum\@xarg>6 \@badlinearg\@xarg 1 \fi
 \setbox\@linechar\hbox{\@linefnt\@getlinechar(\@xarg,\@yyarg)}%
 \ifnum \@yarg >0 \let\@upordown\raise \@clnht\z@
    \else\let\@upordown\lower \@clnht \ht\@linechar\fi
 \@clnwd=\wd\@linechar
 \if@negarg \hskip -\wd\@linechar \def\@tempa{\hskip -2\wd\@linechar}\else
      \let\@tempa\relax \fi
 \@whiledim \@clnwd <\@linelen \do
   {\@upordown\@clnht\copy\@linechar
    \@tempa
    \advance\@clnht \ht\@linechar
    \advance\@clnwd \wd\@linechar}%
 \advance\@clnht -\ht\@linechar
 \advance\@clnwd -\wd\@linechar
 \@tempdima\@linelen\advance\@tempdima -\@clnwd
 \@tempdimb\@tempdima\advance\@tempdimb -\wd\@linechar
 \if@negarg \hskip -\@tempdimb \else \hskip \@tempdimb \fi
 \multiply\@tempdima \@m
 \@tempcnta \@tempdima \@tempdima \wd\@linechar \divide\@tempcnta \@tempdima
 \@tempdima \ht\@linechar \multiply\@tempdima \@tempcnta
 \divide\@tempdima \@m
 \advance\@clnht \@tempdima
 \ifdim \@linelen <\wd\@linechar
    \hskip \wd\@linechar
   \else\@upordown\@clnht\copy\@linechar\fi}
 
 \def\@getlinechar(#1,#2){\@tempcnta#1\relax\multiply\@tempcnta 8
 \advance\@tempcnta -9 \ifnum #2>0 \advance\@tempcnta #2\relax\else
 \advance\@tempcnta -#2\relax\advance\@tempcnta 64 \fi
 \char\@tempcnta}
 
 \def\vector(#1,#2)#3{\@xarg #1\relax \@yarg #2\relax
 \@tempcnta \ifnum\@xarg<0 -\@xarg\else\@xarg\fi
 \ifnum\@tempcnta<5\relax
 \@linelen=#3\unitlength
 \ifnum\@xarg =0 \@vvector 
   \else \ifnum\@yarg =0 \@hvector \else \@svector\fi
 \fi
 \else\@badlinearg\fi}
 
 \def\@svector{\@sline
 \@tempcnta\@yarg \ifnum\@tempcnta <0 \@tempcnta=-\@tempcnta\fi
 \ifnum\@tempcnta <5
   \hskip -\wd\@linechar
   \@upordown\@clnht \hbox{\@linefnt  \if@negarg 
   \@getlarrow(\@xarg,\@yyarg) \else \@getrarrow(\@xarg,\@yyarg) \fi}%
 \else\@badlinearg\fi}
 
 \def\@getlarrow(#1,#2){\ifnum #2 =\z@ \@tempcnta='33\else
 \@tempcnta=#1\relax\multiply\@tempcnta \sixt@@n \advance\@tempcnta
 -9 \@tempcntb=#2\relax\multiply\@tempcntb \tw@
 \ifnum \@tempcntb >0 \advance\@tempcnta \@tempcntb\relax
 \else\advance\@tempcnta -\@tempcntb\advance\@tempcnta 64
 \fi\fi\char\@tempcnta}
 
 \def\@getrarrow(#1,#2){\@tempcntb=#2\relax
 \ifnum\@tempcntb < 0 \@tempcntb=-\@tempcntb\relax\fi
 \ifcase \@tempcntb\relax \@tempcnta='55 \or 
 \ifnum #1<3 \@tempcnta=#1\relax\multiply\@tempcnta
 24 \advance\@tempcnta -6 \else \ifnum #1=3 \@tempcnta=49
 \else\@tempcnta=58 \fi\fi\or 
 \ifnum #1<3 \@tempcnta=#1\relax\multiply\@tempcnta
 24 \advance\@tempcnta -3 \else \@tempcnta=51\fi\or 
 \@tempcnta=#1\relax\multiply\@tempcnta
 \sixt@@n \advance\@tempcnta -\tw@ \else
 \@tempcnta=#1\relax\multiply\@tempcnta
 \sixt@@n \advance\@tempcnta 7 \fi\ifnum #2<0 \advance\@tempcnta 64 \fi
 \char\@tempcnta}
\catcode`\ =10

\dol@g %925
\catcode`\@=\active
\LaTeXfeatures

\newtheorem{theorem}{Theorem}[section]
  \theoremstyle{definition}
\newtheorem{definition }[theorem]{Definition }
\newtheorem{proposition }[theorem]{Proposition }

\newtheorem{corollary}[theorem]{Corollary}

\newtheorem{remark }[theorem]{Remark }

\theoremstyle{remark} \newtheorem{remark}[theorem]{Remark}

\numberwithin{equation}{section}

\begin{document} \title{Automorphisms of categories of free modules and free Lie
algebras.}

\author{R. Lipyanski} \address{Dept. of Mathematics, Ben Gurion
University of the Negev, 84105, Israel}
 \email{lipyansk\@cs.bgu.ac.il}

\author{B. Plotkin} \address{Institute of Mathematics, The Hebrew
University, Jerusalem, 91904, Israel}
\email{plotkin\@macs.biu.ac.il}

\begin{abstract}Let $\Theta^{0}$ be a category of finitly generated
free algebras in the variety of algebras $\Theta$. Solutions to
problems in algebraic geometry over $\Theta$ are often determined
by the structure of the group of automorphisms $Aut\,\Theta^{0}$
of category $\Theta^{0}$. Here we consider two varieties $\Theta$:
noetherian modules and Lie algebras. We show that every
automorphism in $Aut\,\Theta^{0}$, where $\Theta$ is the variety
of modules over noetherian rings, is semi-inner. A similar result
for the variety of Lie algebras over a infinite field has been
recently obtained in \cite{Mash1}. Here we present a different
approach allowing us to shorten the proof in the Lie case.
  \end{abstract}
 \maketitle

 \section{Introduction}
 \subsection{ Automorphisms of category. }

 Let $\bf {\Phi}$ be a small category. Denote by $End\,\bf {\Phi}$ the semigroup of
 endofunctors  of the category $\bf{\Phi}$ and by $Aut\,\bf {\Phi}$ the subgroup of
 invertible elements of $End\,\bf{\Phi}$.

Recall that functor isomorphism
 $s : \varphi_{1}\longrightarrow \varphi_{2}$,\;\;$\varphi_{1}, \varphi_{2}\in
 End\,\bf {\Phi}$, is a collection of isomorphisms
  $s_{A}: \varphi_{1}(A)\longrightarrow \varphi_{2}(A)$ defined for
  all $A\in Ob\,\bf {\Phi}$ and
  such that for every $\nu:A\longrightarrow B$,\; $\nu\in Mor\,\bf
   {\Phi}$, holds
   $$
     s_{B}\cdot\varphi_{1}(\nu)= \varphi_{2}(\nu)\cdot s_{A},
   $$
  i.e., the following diagram is commutative

 $$\CD \varphi_{1}(A) @>s_{A} >> \varphi_{2}(A)\\ @V\varphi_{1}(\nu)
  VV @VV\varphi_{2}(\nu) V\\\varphi_{1}(B)@>s_{B}>> \varphi_{2}(B)\endCD$$

  A relation $\simeq$ of the isomorphism on functors of a category $\bf {\Phi}$
  is a congruence on the semigroup $End\,\bf {\Phi}$. Denote
  $ End^{\circ}\,\Phi=End\,\Phi\,/\simeq $.
  We have a natural homomorphism $\delta:End\,\Phi\longrightarrow End^{\circ}\,\bf {\Phi}$
  which induces a homomorphism $Aut\,{\bf \Phi} \longrightarrow
  Aut^{\circ}\,\bf \Phi$, were $Aut^{\circ}\,\bf{\Phi}$ is the group of
  invertible elements of the semigroup $End^{\circ}\,\bf{\Phi}$.
  \begin{definition }
  An automorphism $\varphi\in Aut^{\circ}\,\bf{\Phi}$ of the category
   $\bf{\Phi}$ is isomorphic to the identical one is called inner.
  \end{definition }
  That means an automorphism $\varphi:\bf{\Phi}\longrightarrow \bf\Phi$
  is inner if there exists a function $s$ which chooses for each
  object $A\in Ob\,\bf\Phi$ an isomorphism $s_{A}: A\longrightarrow \varphi(A)$
  such that for every $\nu: A\longrightarrow B$ holds

  $$\CD A @>s_{A} >> \varphi(A)\\ @V\nu VV @VV\varphi (\nu) V\\
B @>s_{B} >> \varphi(B)\endCD$$
 Here $\varphi(\nu)=s_{B}\cdot\nu\cdot s^{-1}_{A}$.
 It is clear that the kernel $Ker\,\nu$ of the homomorphism $\delta$ consists all
  inner automorphisms of the category $\bf\Phi$.

  Denote by $Int\,\Phi$ the group of inner automorphisms of the
  category $\Phi$.
    \begin{definition }
    Group $Out\,{\bf\Phi}=Aut\,{\bf\Phi}/Int\,{\bf\Phi}$ is called the
    group of outer automorphisms of the category $\Phi$.
     \end{definition }
 \begin{definition }
An automorphism $\varphi$ of the category $\bf\Phi$ is called a
special automorphism of this category if for each object $A\in
\bf\Phi$ we have $A\simeq \varphi(A)$.
 \end{definition }
  \begin{definition }
  The category is called special if each of its automorphisms is
  special.
    \end{definition }
     \begin{definition }
  The category is called perfect if each of its automorphisms is
 inner.
    \end{definition }
 It is easy to prove
 \begin{proposition }\label{pr0}
 If $\varphi$ is a special automorphism of the category $\bf\Phi$
 then it can be represented in the following form: $\varphi
 =\varphi_{0}\cdot\varphi_{1}$, where $\varphi_{0}$ is an inner
 automorphism of the category $\bf{\Phi}$ and $\varphi_{1}$ is an
 automorphism which changes no object of $\bf\Phi$.
  \end{proposition }

 Bellow we consider the category $\Theta^0$ of free algebras $W=W(X)$ generated by
finite sets $X$, where each $X$  belongs to the universe $X^{0}$.
In this case the category $\Theta^{0}$ is a small category.

It has been proven \cite{Zhit}  that the homomorphism $\delta :
Aut\,\Theta^{0}\longrightarrow Aut^{0}\,\Theta^{0}$ is a
surjective.

Varieties of three following types are important for our
consideration below.
 \begin{definition }\label{p1}
  A variety $\Theta$ is notherian if each object of the category
  $\Theta^{0}$ is notherian with respect to its congruences.
  \end{definition }
       \begin{definition } \label{dp1}
 A variety $\Theta$ is hopfian if  each object of the category
 $\Theta^{0}$ is hopfian, i.e. every surjective endomorphism $\mu : W\longrightarrow W$ is
an automorphism.
     \end{definition }
      \begin{definition }\label{dp2}
 A variety $\Theta$ is regular if the condition
$W(X)\simeq W(Y)$ implies $|X|\simeq |Y|$.
    \end{definition }
The following implications are true: $\ref{p1}\Longrightarrow
\ref{dp1}\Longrightarrow \ref{dp2}$.

We recall several results related to the category $\Theta^{0}$.
 \begin{theorem}\cite{Mash2} \label{t1}
 Let $W_{0}=W(x_{0})$ be a cyclic algebra generated by $x_{0}$,
 $\Theta$ a hopfian variety and $\varphi$ an automorphism of
  $\Theta^{0}$. Then from $\varphi(W_{0})\simeq W_{0}$ follows that
 $\varphi$ is a special automorphism.
 \end{theorem}
 Let $\Theta$ be a variety generated by a free algebra
 $W^{0}=W(x_{1},...,x_{k})$ and $\nu_{0}: W^{0}\longrightarrow W_{0}$
  be a homomorphism determined by the equalities
 $\nu_{o}(x_{i})=x_{0},\,i=1,...,k$.
   \begin{theorem}\cite{Mash2}\label{t2}
     Let $\Theta$ be a hopfian variety and $\varphi$ an
     automorphism of the category $\Theta^{0}$ not changing its
     objects and acting on $End\,W^{0}$ as an identical
     automorphism. In addition let  $\varphi(\nu_{0})=\nu_{0}$.
     Then $\varphi$ is an inner automorphism.
   \end{theorem}
 Note that any automorphism $\varphi$ of the category $\Theta^{0}$
 which does not changes objects of $\Theta^{0}$ induces an
 automorphism $\varphi_{W}$ of the semigroup $End\,W$ for every
 $W=W(X),$ $|X| <\infty$. Furthermore, $\varphi$ induces a
 substitution on the set $Hom\,(W_{0},W)$. Denote this
 substitution by $\sigma_{X}$.
  \begin{theorem}\cite{Mash2}\label{t3}
     Let $\nu : W(X)\longrightarrow W(Y)$ be a homomorphism. Then
     $\varphi(\nu)=\sigma_{Y}\cdot\nu \cdot\sigma^{-1}_{X}$.
   \end{theorem}
 \subsection{Geometries over algebras}
 Recall the main definitions from the universal algebraic geometry
 \cite{Plot}. Let $\Theta$ be a variety of algebras and $W=W(X)\in \Theta$
 be a free algebra from $\Theta$ over  a finite set of generators
$X,\,|X|=n $. The set $Hom(W,\,G),\,G\in\Theta$ can be treated as
an affine space whose points are homomorphisms. Every solution of
an equation $w=w^{\prime}$ in $W$ is a homomorphism $\mu :
W\longrightarrow G$ such that $w^{\mu}=w^{\prime \mu}$. The set
$A\subseteq Hom(W,G)$ is called an algebraic set over the algebra
$G$ if there exists a system of equations $T$ (or a binary
relation $T$ on $W$) such that $T^{\prime}=A$, where
$T^{\prime}=\{\mu\mid T\subseteq Ker\,\mu \}$. It is possible to
define also the relation $A^{\prime}$ on $W$:
$A^{\prime}=T=\bigcap\limits_{\mu\in A} Ker\,\mu$ and if $A$ is
an algebraic set then $A^{\prime\prime}=A$.

 Denote by $K_{\Theta}(G)$ the category of pairs $(A,X)$, where
 $A$ is an algebraic set over $G$, i.e., $A\subseteq Hom(W(X),G)$.
 Morphisms in this category are mappings
 $$
  (\overline{s},\,s):(A,X)\longrightarrow (B,Y),
 $$
 where $s:W(Y)\longrightarrow W(X)$ is a mapping in category
 $\Theta^{o}$ of free algebras from ${\Theta}$ and
 $\overline{s}:A\longrightarrow B$ is such a mapping from $A$ to $B$
 induced by $s$ that the following diagram  is  commutative:

   $$\CD W(Y) @>s >> W(Y)\\ @V\mu_{Y}VV @VV\mu_{X} V\\
W(Y)/B^{\prime} @>\overline{s}>> W(X)/A^{\prime} \endCD$$

 Here $\mu_{X}$ and  $\mu_{Y}$ are the natural homomorphisms.
 The category $K_{\Theta}(G)$ is a geometric invariant of $G$.
 Algebras $G_{1}$ and $G_{2}$ are categorically equivalent
 if the categories $K_{\Theta}(G_{1})$ and $K_{\Theta}(G_{2})$ are
 isomorphic. If algebras $G_{1}$ and $G_{2}$ are categorically equivalent
 then these algebras define the same geometry.
  Algebras $G_{1}$ and $G_{2}$ are geometrically equivalent if

  $$
  T^{\prime\prime}_{G_{1}}= T^{\prime\prime}_{G_{2}}
    $$
  holds for all finite sets $X$ and for all relations $T$ in $W=W(X),\,W\in
  \Theta$.

   In some cases the categorical equivalence of algebras $G_{1}$ and $G_{2}$
   coincides with the geometric equivalence of these algebras.
   Whether these two notions are or are not equivalent in $\Theta$, depends on the
   structure of the group $Aut\,\Theta^{0}$ \cite{Plot}.
 The group $Aut\,\Theta^{0}$ is known in the
 following cases: the variety of all groups, the variety of $F$-groups, where
 $F$ is a free group of constants, the variety of all semigroups, the variety of
 commutative associative algebras with unit element over an
 infinite field, the variety of all Lie algebras over an infinite field
 \cite{Ber}, \cite{Mash1}, \cite{Mash2},  \cite{Plot}.

  Here we give a description of the group $Aut\,\Theta^{0}$ for the variety
 of all $K$-modules, where $K$ is a noeterian ring, and  present a
 different approach allowing to simplify the proof \cite{Mash1} of the corresponding
 theorem for the variety of Lie algebras over an infinite field.
 \section{Modules}
   \subsection{Preliminary remarks}
  Let $K$ be an arbitrary ring with unit. Denote by
  $\bf C=Mod-K$ the category of all left $K$-modules. Every free
  module is a direct sum of cyclic modules over $K$. The ring $K$
  is noetherian if it is noetherian as a $K$-module. Thus the
  variety $\bf C$ is noetherian iff the ring $K$ is
  noetherian. In this case the variety $\bf C$ is also hopfian and
  regular.

   We wish to study the semigroup of endomorphisms of a cyclic module over $K$.
   Let $F=Kx_{0}$ be a free cyclic $K$-module generated  by $x_{0}\in X$ and
   $\nu_{\alpha} : Kx_{0}\longrightarrow Kx_{0}$ be an
  endomorphism of this module such that $\nu_{\alpha}(x_{0})=\alpha x_{0},\,\alpha\in
  K$. Since
  $(\nu_{\alpha}\cdot\nu_{\beta})(x_{0})=\nu_{\beta\cdot\alpha}(x_{0})$, the
  semigroup $End\,Kx_{0}$ is antiisomorphic to the multiplicative
  semigroup $K^{\times}$ of the ring $K$. We will show first that
  there exists isomorphism $ Aut\,(End\,Kx_{o})\simeq Aut\,K^{\times}$
  between groups of automorphisms of the semigroups $End\,Kx_{o}$ and $K^{\times}$ respectively.

  Let $\sigma$ be an automorphism of the semigroup $K^{\times}$,
   i.e.,
  $(\beta\cdot\alpha)^{\sigma}=\beta^{\sigma}\cdot\alpha^{\sigma}$, $\alpha,\beta\in K^{\times}$.
 Define a mapping $\sigma^{*}$ on $End\,Kx_{o}$ in the following
 way: $\nu_{\alpha}^{\sigma^{*}}=v_{\alpha^{\sigma}}$.
  We may write
  $$
      (v_{\alpha}\cdot
      v_{\beta})^{\sigma^{*}}=v_{\beta\alpha}^{\sigma^{*}}=v_{(\beta\alpha)^{\sigma}}
      =v_{\beta^{\sigma}\alpha^{\sigma}}
      = v_{\alpha}^{\sigma^{*}}\cdot
      v_{\beta}^{\sigma^{*}},
  $$
  i.e., $\sigma^{*}\in Aut\,(End\,Kx_{0})$.

  Next we assume that $\sigma^{*}\in Aut\,(End\,Kx_{0})$, i.e.,
  ${(v_{\alpha}\cdot v_{\beta})}^{\sigma^{*}}=
  v_{\alpha}^{\sigma^{*}}\cdot v_{\beta}^{\sigma^{*}}$, for all
  $v_{\alpha},v_{\beta}\in End\,Kx_{0}$.
   Then there exists a mapping $\sigma: K^{\times}\longrightarrow K^{\times}$
   such that
  $\nu_{\alpha}^{\sigma^{*}}=v_{\alpha^{\sigma}}$. We will show
  that $\sigma\in Aut\,K^{\times}$.

  We have
  $$
   v_{\alpha}^{\sigma^{*}}
   v_{\beta}^{\sigma^{*}}=v_{\alpha^{\sigma}}v_{\beta^{\sigma}}=
   v_{\beta^{\sigma}\alpha^{\sigma}}.
    $$
 On the other hand
 $$
   {(v_{\alpha} v_{\beta})}^{\sigma^{*}}=v_{\beta\alpha}^{\sigma^{*}}=
   v_{{(\beta\alpha)}^{\sigma}}.
   $$
Thus
$(\beta\cdot\alpha)^{\sigma}=\beta^{\sigma}\cdot\alpha^{\sigma}$
and $\sigma\in Aut\,K^{\times}$. It is clear that
$({\sigma_{1}\cdot\sigma_{2}})^{*}=
{\sigma_{1}}^{*}\cdot{\sigma_{2}}^{*}$.

 \begin{proposition }\label{p3}
 Let $\bf C^{0}=(Mod-K)^{0}$ be a category of finite generated free modules
 over noetherian ring $K$. If $\varphi$ is an automorphism of
 the category $\bf C^{0}$ then $\varphi_{Kx_{0}}=\sigma$ determines
 an automorphism of the ring $K$.
 \end{proposition }
 \begin{proof}
  Let $A=Kx_{1}\oplus
Kx_{2},\,x_{1},x_{2}\in X$ be a 2-generated free $K$-module and
$s^{\alpha}$ be  its automorphism of the following form: $$
\begin{array}{l}
 s^{\alpha}(x_{1})=x_{1}+\alpha\cdot x_{2},\quad\alpha\in K \\
 s^{\alpha}(x_{2})=x_{2}
  \end{array}
 $$
 Denote by the same symbol the matrix of the automorphism
 $s^{\alpha}$:
  $$
   s^{\alpha}=\left(\begin{array}{cc}1& \alpha\\ 0& 1
    \end{array}   \right)
  $$
 Consider the following commutative diagram
 $$
    \CD
Kx_{0} @>\varepsilon_1 >> Kx_{1}\oplus Kx_{2}\\
@. @/SE/\delta_1 // @VV\pi_1 V\\
@. Kx_1\\
\endCD $$ \noindent
 where $\delta_{1}=\pi_{1}\cdot\varepsilon_{1},\;\varepsilon_{1}(x_{0})=x_{1}$
 and $\pi_{1}$ is the natural projection $Kx_{1}\oplus Kx_{2}$ on
 $Kx_{1}$, i.e., $\pi_{1}(x_{1})=x_{1},\;\pi_{1}(x_{2})=0$.
 Applying the automorphism $\varphi$ to the last diagram, we obtain
  $$
\CD
Kx_{0} @>\varphi(\varepsilon_1) >> Kx_{1}\oplus Kx_{2}\\
@. @/SE/\varphi(\delta_1 )// @VV\varphi(\pi_1)V\\
@. Kx_1\\
\endCD $$ \noindent
 Let $\varphi(\varepsilon_{1})(x_{0})=y_{1}$ and $\varphi(\pi_{1})(y_{1})
 =\varphi(\delta_{1})(x_{0})=y_{1}^{0}\in Kx_{1}$.

 Now, we consider another commutative diagram:
 $$
    \CD
Kx_{0} @>\varepsilon_2 >> Kx_{1}\oplus Kx_{2}\\
@. @/SE/\delta_2 // @VV\pi_2 V\\
@. Kx_2\\
\endCD $$
where
$\delta_{2}=\pi_{2}\cdot\varepsilon_{2},\;\varepsilon_{2}(x_{0})=x_{2}$
 and $\pi_{2}$ is the natural projection $Kx_{1}\oplus Kx_{2}$ on
 $Kx_{2}$, i.e., $\pi_{2}(x_{2})=x_{2},\;\pi_{2}(x_{1})=0$.
 Applying $\varphi$ to the last diagram,
 we obtain
 $$
  \varphi(\varepsilon_{2})(x_{0})=y_{2},\; \varphi(\pi_{2})(y_{2})
 =\varphi(\delta_{2})(x_{0})=y_{2}^{0}\in Kx_{2}.
   $$
   Since the ring $K$ is noetherian, we receive
   $$
    Kx_{1}\oplus Kx_{2}=Ky_{1}\oplus Ky_{2}
   $$
Now, consider the diagram
   $$
    \CD
Kx_{0} @>\varepsilon_1 >> Kx_{1}\oplus Kx_{2}\\
@. @/SE/\delta // @VV\pi_2 V\\
@. Kx_2\\
\endCD $$ \noindent
 where $\pi_{2}\cdot\varepsilon_{1}=\delta$. It is clear that
 $\pi_{1}\cdot\varepsilon_{2}=0$ and so
 $$
   \varphi(\pi_{1})\cdot\varphi(\varepsilon_{2})=0 \quad\mbox{and}\quad
   \varphi(\pi_{2})\cdot\varphi(\varepsilon_{1})=0.
   $$
 Thus
     $$
 \varphi(\pi_{1})\cdot\varphi(\varepsilon_{2})(x_{0})= \varphi(\pi_{1})(y_{2})=0
 \quad\mbox{and}\quad
 \varphi(\pi_{2})\cdot\varphi(\varepsilon_{1})(x_{0})=
 \varphi(\pi_{2})(y_{1})=0.
     $$
  We obtained the following relations:
   \begin{equation}\label{forr1}
  \begin{array}{l}
 \varphi(\pi_{1})(y_{2})=0,\; \varphi(\pi_{1})(y_{1})=y_{1}^{0}\in Kx_{1},\\
    \varphi(\pi_{2})(y_{1})=0,\; \varphi(\pi_{2})(y_{2})=y_{2}^{0}\in Kx_{2}.
  \end{array}
    \end{equation}

 By definition of $s^{\alpha}$
  $$
  \begin{array}{l}
  s^{\alpha}(x_{1})=s^{\alpha}\varepsilon_{1}(x_{0})=x_{1}+\alpha\cdot x_{1},\\
      \pi_{1}s^{\alpha}\varepsilon_{1}(x_{0})=x_{1}=\delta_{1}(x_{0}).
       \end{array}
 $$
Thus $\pi_{1}s^{\alpha}\varepsilon_{1}=\delta_{1}$ and so
$\varphi(\pi_{1})\varphi(s^{\alpha})\varphi(\varepsilon_{1})=\varphi(\delta_{1})$.
Applying the last equality to $x_{0}$, we get
  \begin{equation}\label{forr2}
 \varphi(\pi_{1})\varphi(s^{\alpha})(y_{1})=\varphi(\delta_{1})(x_{0})
 =y^{0}_{1}=\varphi(\pi_{1})(y_{1}).
  \end{equation}
  Further we have
  $$
   \pi_{2}s^{\alpha}\varepsilon_{1}(x_{0})=\alpha\cdot x_{2}=\alpha\cdot\delta_{2}(x_{0})
     =\delta_{2}(\alpha\cdot x_{0})=\delta_{2}\nu_{\alpha}(x_{0}).
    $$
   Therefore $\pi_{2}s^{\alpha}\varepsilon_{1}= \delta_{2}\nu_{\alpha}$, hence
     \begin{equation} \label{forr3}
      \varphi(\pi_{2})\varphi(s^{\alpha})\varphi(\varepsilon_{1})=
      \varphi(\delta_{2})\varphi(\nu_{\alpha})=
      \varphi(\delta_{2})\nu_{\alpha^{\sigma}}.
     \end{equation}
   Here $\sigma$ is an automorphism of the semigroup $K^{\times}$
   corresponding to $\varphi$.
  Now, applying equalities (\ref{forr3}) to the element $x_{0}$, we
 obtain
     \begin{equation} \label{forr4}
    \begin{array}{l}
      \varphi(\pi_{2})\varphi(s^{\alpha})(y_{1})=
     \varphi(\delta_{2})\nu_{\alpha^{\sigma}}(x_{0})=
        \\= \varphi(\delta_{2})(\alpha^{\sigma}x_{0})=
      \alpha^{\sigma}  \varphi(\delta_{2})(x_{0})= \alpha^{\sigma}
      y^{0}_{2}.
      \end{array}
      \end{equation}
   We come to
      $$
     \varphi(\pi_{2})\varphi(s^{\alpha})(y_{1})= \alpha^{\sigma}
      y^{0}_{2}= \alpha^{\sigma}\varphi(\pi_{2})(y_{2}).
       $$
  Further we get
      $$
  \begin{array}{l}
      s^{\alpha}(x_{2})=x_{2},\; s^{\alpha}\varepsilon_{2}(x_{0})=\varepsilon_{2}(x_{0})\\
     s^{\alpha}\varepsilon_{2}= \varepsilon_{2},\;
     \varphi(s^{\alpha})\varphi(\varepsilon_{2})=
     \varphi(\varepsilon_{2}).
  \end{array}
      $$
    Applying the last equality to $x_{0}$, we obtain
    \begin{equation} \label{forr5}
    \varphi(s^{\alpha})(y_{2})=y_{2}.
       \end{equation}
   Now we need to represent the element $\varphi(s^{\alpha})(y_{1})$
  in the basis $y_{1},\;y_{2}$. Let
  $\varphi(s^{\alpha})(y_{1})=\lambda_{1}\cdot y_{1}+\lambda_{2}\cdot
  y_{2}$. By the equalities (\ref{forr1}), (\ref{forr2}),
  (\ref{forr4}), (\ref{forr5}) we calculate $\lambda_{1}$
  $$
   \begin{array}{l}
    \varphi(\pi_{1})\varphi(s^{\alpha})(y_{1})=y_{1}^{0}=\\=
     \lambda_{1}\varphi(\pi_{1})
     (y_{1})+\lambda_{2}\varphi(\pi_{2})(y_{2})=
     \lambda_{1}{y_{1}^{0}}, \;\mbox{i.e.}\;\,\lambda_{1}=1,
     \end{array}
    $$
  and $\lambda_{2}$
  $$
   \begin{array}{l}
    \varphi(\pi_{2})\varphi(s^{\alpha})(y_{1})= \varphi(\pi_{2})(y_{1}
    + \lambda_{2}
    y_{2})=\\=\lambda_{2}\varphi(\pi_{2})(y_{2})=\lambda_{2}y_{2}^{0}=\alpha^{\sigma}
    y_{2}^{0}, \;\mbox{i.e.}\;\,\lambda_{2}=\alpha^{\sigma}.
     \end{array}
   $$
   Finally, we obtain
     $$
   \begin{array}{l}
      \varphi(s^{\alpha})(y_{1})=y_{1}+ \alpha^{\sigma}y_{2}     \\
      \varphi(s^{\alpha})(y_{2})=y_{2}.
     \end{array}
   $$

We showed that the matrix of $\varphi(s^{\alpha})$ relative to the
basis $y_{1},y_{2}$ has the same triangular form as $s^{\alpha}$:

  $$
   s^{{\alpha}^{\sigma}}=\left(\begin{array}{cc}1& \alpha^{\sigma} \\ 0& 1
    \end{array}   \right).
  $$
  We have $s^{\alpha}s^{\beta}=s^{\alpha+\beta}$. Since
  $$
     \varphi(s^{\alpha})\varphi(s^{\beta})=\varphi(s^{\alpha+\beta}),
  $$
 we must have
   $$
   s^{\alpha^{\sigma}}s^{\beta^{\sigma}}=s^{\alpha^{\sigma}+\beta^{\sigma}}=
      s^{{(\alpha+\beta)}^{\sigma}}.
     $$
 Thus ${(\alpha+\beta)}^{\sigma}=\alpha^{\sigma}+\beta^{\sigma}$ and $\varphi_{Kx_{0}}=\sigma$
  is an automorphism of the ring $K$.

 \end{proof}
 \subsection{Semi-inner automorphism of the category of free modules over a ring}
Let us define the notion of a semi-inner automorphism of the
category $\bf C^{0}=(Mod-K)^{0}$.
  Consider the category of modules with semimorphisms
 (semi-linear maps).
   \begin{definition }
   A semimorphism of the category $\bf Mod-K$ is a pair $(\sigma,
   s)$, where $\sigma\in Aut\,K$ and $s:A\longrightarrow B$ an
   addition-compatible map such that
   $s(\alpha\cdot a)=\alpha^{\sigma}s(a),\;a\in A,\;\alpha\in K$.
     \end{definition }

  \begin{definition }
  An automorphism $\varphi$ of the category $\bf(Mod-K)^{0}$ is
  called semi-inner if there exists  $\sigma\in Aut\,K$ and a
  collection
  $(\sigma,s)=\{(\sigma,s_{X})|s_{X}:KX\longrightarrow KX,\, \sigma\in Aut\,K\}$
   of semi-isomorphisms such that $(\sigma, s) :1_{\bf C^{0}}\longrightarrow \varphi$ is
  semi-isomorphism of functors, i.e.,  for every object $KX\in Ob\,\bf C^{0}$
 there exists semi-isomorphism $(\sigma, s_{X}): KX\longrightarrow
 \varphi(KX)$ and for each $\nu: KX\longrightarrow KY$ the
 following diagram is commutative:
    $$\CD KX @>s_{X} >> \varphi(KX)\\ @V\nu VV @VV\varphi (\nu) V\\
KY @>s_{Y} >> \varphi(KY)\endCD$$
  \end{definition }
   Let $\sigma$ be an automorphism of the ring $K$ and $KX=Kx_{1}\oplus...\oplus Kx_{n}$
 be a free module over a ring $K$. Define a mapping $\sigma_{X}:KX\longrightarrow KX$
 such that if $u=\sum\limits_{i=1}^{n}\lambda_{i}x_{i}$ then
 $\sigma_{X}(u)=\sum\limits_{i=1}^{n}\lambda_{i}^{\sigma}x_{i}$.
 It is evident that
 $$
  \sigma_{X}(u+v)=\sigma_{X}(u)+\sigma_{X}(v)\,\;\mbox{and}\;\,
     \sigma_{X}(\lambda u)= \lambda^{\sigma}\sigma_{X}(u),
  $$
for each $u,v\in KX$ and $\lambda\in K$.

  Now, consider a mapping
  $$
  \varphi(\nu)=\sigma_{Y}\nu\,\sigma_{X}^{-1} :KX\longrightarrow KY
  .
  $$
We show that this mapping is a homomorphism of modules. Clearly,
it preserves the additive structure. Taking $\lambda u\in
KX,\,\lambda\in K,\,u\in KX$, we obtain $$
  \begin{array}{l}
   \sigma_{Y}\nu\,\sigma_{X}^{-1}(\lambda u)=
   \sigma_{Y}\nu(\lambda^{\sigma^{-1}}\sigma^{-1}_{X}(u))=\\
   = \sigma_{Y}\lambda^{\sigma^{-1}}\nu\sigma^{-1}_{X}(u)=
   \lambda \sigma_{Y}\nu\sigma^{-1}_{X}(u),\\
     \end{array}
 $$
 i.e., $\sigma_{Y}\nu\,\sigma_{X}^{-1}(\lambda u)
 =\lambda \sigma_{Y}\nu \sigma^{-1}_{X}(u)$. Therefore
 $\sigma_{Y}\nu\,\sigma_{X}^{-1}$ is a module homomorphism.

 Denote by $\hat{\sigma}$ an automorphism of the category of
 free $K$-modules of the following form: \par\medskip
 1. \,$\hat{\sigma}(KX)=KX, \;KX\in Ob\;\bf C^{0}$,\, i.e.,
  $\hat{\sigma}$ does not change objects
  of the category of the free $K$-modules. \par\medskip
 2. \,$\hat{\sigma}(\nu)=\sigma_{Y}\nu\,\sigma_{X}^{-1}:KX\longrightarrow KY $
  for all $\nu:KX\longrightarrow KY;\; KX,\,KY\in Ob\;\bf C^{0}$.

  It is clear that $\hat{\sigma}$ is
  a semi-inner automorphism of the category $\bf C^{0}$.

 Note that the set $S$ of the semi-inner automorphisms of the category
  $\bf C^{0}$ is a subgroup in $Aut\,\bf C^{0}$. The subgroup $S$ contains the invariant
 subgroup of the inner automorphisms of the category $\bf C^{0}$.

 Now, we shall describe the automorphisms of the category $\bf(Mod-K)^{0}$
over a noetherian  ring $K$.
 \begin{theorem}\label{main1}
  Let $K$ be a noetherian  ring. Every automorphism of the
  category $\bf C^{0}=(Mod-K)^{0}$ is semi-inner.
  \end{theorem}
 \begin{proof}
 Suppose that $\varphi\in Aut\,\bf{C}^{0}$. We shall show that $\varphi(Kx_{0})\simeq
 Kx_{0}$. Recall that $ Aut\,(End\,Kx_{o})\simeq Aut\,K^{\times}$.
 It is clear that $ Aut\,(End\,KX)\not\simeq Aut\,K^{\times}$ if
 $|X|\neq 1$. Therefore $\varphi(KX)\not\simeq Kx_{0}$ if $|X|\neq 1$
 and $\varphi(Kx_{0})\simeq Kx_{0}$.

  Since $K$ is a noetherian ring, the variety  $\bf C$ is hopfian.
  By Theorem \ref{t1} every automorphism of the
 category $\bf C^{0}$ is special. By Proposition \ref{pr0} the automorphism
 $\varphi$ can be represented as
 $\varphi=\varphi_{1}\cdot\varphi_{2}$, where $\varphi_{1}$ is an inner
 automorphism of $\bf C^{0}$ and $\varphi_{2}$ is an automorphism which does not change
  the objects of this category. Now, without loss of
  generality, we can assume that $\varphi$ itself does not change the objects of
  the category $\bf C^{0}$.

  Take $\varphi_{Kx_{0}}=\sigma$. Denote by the same letter
  the automorphism of the ring $K$ corresponding to $\sigma$
  (see Proposition \ref{p3}). Let $\hat{\sigma}$ be a semi-inner automorphism of
   the category $\bf C^{0}$ corresponding to $\sigma$. We shall show that
   $\hat{\sigma}$ and $\sigma$ act in the same way on $End\,Kx_{0}$.

  Let $\nu_{\alpha}: Kx_{0}\longrightarrow Kx_{0}$ be a an
endomorphism of the cyclic module $Kx_{0}$. Consider the
endomorphism $\sigma_{x_{0}}\nu_{\alpha}\sigma_{x_{0}}^{-1}$. Here
for the convenience we denote by $x_{0}$ the subset $\{x_{0}\}$.
We have $$
 \sigma_{x_{0}}\nu_{\alpha}\sigma_{x_{0}}^{-1}(x_{0})
 =\sigma_{x_{0}}\nu_{\alpha}(x_{0})= \sigma_{x_{0}}(\alpha
 x_{0})=\alpha^{\sigma}x_{0}=\nu_{{\alpha}^{\sigma}}(x_{0}).
$$ Therefore
 $$
   \sigma_{x_{0}}\nu_{\alpha}\sigma_{x_{0}}^{-1}=
   \nu_{{\alpha}^{\sigma}}=\nu_{\alpha}^{\sigma}
 $$
Thus $\hat{\sigma}_{Kx_{0}}=\sigma$ and
$\varphi_{1}=\hat{\sigma}^{-1}\varphi$ is an automorphism of $\bf
C^{0}$ acting identically on the semigroup $End\,(Kx_{0})$.

   Recall that the variety $\bf C$ generated by a cyclic module
   $Kx_{0}$. Consider a homomorphism $\nu_{0}:Kx_{0}\longrightarrow
  Kx_{0}$ such that $\nu_{0}(x_{0})=x_{0}$. It is clear that
   $\varphi_{1}(\nu_{0})=\nu_{0}$. Hence all conditions of
   Theorem \ref{t2} hold, and the automorphism  $\varphi_{1}=\hat{\sigma}^{-1}\varphi$
   is an inner automorphism of the category $\bf C^{0}$. Thus
    $ \varphi=\hat{\sigma}\varphi_{1}$ is semi-inner. The proof is complete.
    \end{proof}
   \begin{corollary}
    If a ring $K$ has non-trivial automorphisms, the
    variety $Mod-K$ is perfect, i.e., all automorphisms of $\bf (Mod-K)^{0}$
    are inner.
     \end{corollary}

    \begin{corollary}
    The variety of abelian groups is perfect.
     \end{corollary}
   Let now $K$ be an arbitrary ring. Recall that the group $Out\,K$ of outer
  automorphisms of the ring $K$ is a factor group
  $Aut\,K/Int\,K$, were $Aut\,K$ and $Int\,K$ are the groups of all
  automorphisms and the group of inner automorphisms respectively
  of the ring
  $K$.
    \begin{theorem}
  If automorphisms of the variety $\bf (Mod-K)^{0}$ are semi-inner then the group of
  outer automorphisms of the category ${\bf C}=\bf (Mod-K)^{0}$ is isomorphic to the
  group of outer automorphisms of the  ring $K$.
 \end{theorem}
 \begin{proof}
    We should prove the isomorphism
  $$
  Aut\,{\bf C}/Int\,{\bf C}\simeq Out\,K
      $$

Denote by $Sp\,\,{\bf C}$ the group of the all special
automorphisms of the category ${\bf C}$ and by $St\,\,{\bf C}$ the
group of the all automorphisms which do not change objects of
${\bf C}$.
 We have
 $$
 \begin{array}{l}
    Sp\,\,{\bf C}= St\,{\bf C}\cdot Int\,{\bf C}\\
   Sp\,\,{\bf C}/Int\,{\bf C}\simeq St\,{\bf C}/Int\,{\bf C}\cap
    St\,{\bf C}.
     \end{array}
  $$
 As a consequence of these formulas we get
    $$
 \begin{array}{l}
    Sp\,\,{\bf C}= Aut\,{\bf C}\\
    Aut\,{\bf C}/Int\,{\bf C}\simeq St\,{\bf C}/Int\,{\bf C}\cap
    St\,{\bf C}.
     \end{array}
  $$
 Now, it is sufficient establish an isomorphism
    $$
   St\,{\bf C}/Int\,{\bf C}\cap St\,{\bf C}\simeq Out\,K .
      $$
   Let us construct a group homomorphism
   $$
       \pi: St\,{\bf C}\longrightarrow Aut\,K .
     $$
 Any automorphism $\varphi\in St\,{\bf C}$ induces an automorphism $\varphi_{Kx_{0}}$
  of the semigroup $End\,Kx_{0} $. Since $ Aut\,(End\,Kx_{o})\simeq
  Aut\,K^{\times}$, there exists $\sigma\in Aut\,K^{\times}$
  corresponding to $\varphi_{Kx_{0}}$. Put $\pi(\varphi)=\sigma$.

Since $\pi(\hat{\sigma})=\sigma$ for all $\sigma\in Aut\,K $,
homomorphism $\pi$ is surjective. Now consider the natural
homomorphism $\psi$
   $$
  \psi: Aut\,K\longrightarrow Aut\,K/ Int\,K.
   $$
Denote by $\overline{\pi}= \psi\cdot\pi$. We have to describe the
kernel $Ker\,\overline{\pi}$. By definition $\varphi\in
Ker\,\overline{\pi}$ iff the automorphism $\sigma$ of the ring
$K$ corresponding to $\varphi$ is an inner automorphism of this
ring. Equivalently, $\varphi_{Kx_{0}}$ is an inner automorphism of
the semigroup $End\, Kx_{0}$. By Theorem \ref{t3} we have that
$\varphi\in Int\, {\bf C}\cap St\,{\bf C}$. Finally, we get
 $$
  Ker\,\overline{\pi}=Int\,{\bf C}\cap St\,{\bf C}.
 $$
This completes the proof of the theorem.
  \begin{corollary}
  If $K$ is a noetherian ring then the group of the outer automorphisms of the category
  $\bf (Mod-K)^{0}$, is isomorphic to the group of the outer automorphisms of the  ring $K$.
   \end{corollary}
  \end{proof}
  \section{Lie algebras}
    \subsection{Automorphisms of semigroup of endomorphisms of 2-generated Lie algebras}
We consider the variety $\bf\Theta=Lie-K$ of all Lie algebras over
field $K$.
   \begin{definition } \label{dl1}
  Let $W=W(X)$ be a free Lie algebra over a field $K$ and $\delta\in Aut\,K$. Denote a
  mapping $\delta_{W}$ of the algebra $W$ to $W$ in the following
  way: \par\medskip
  1. \,$ \delta_{W}([w_{1},w_{2}])=[\delta_{W}(w_{1}),\delta_{W}(w_{2})],\,
  w_{1},w_{2}\in W,$ \par\medskip
  2. \,$\delta_{W}(\lambda \cdot w)=\lambda^{\delta}\cdot\delta_{W}(w),\,w\in W,
  \,\lambda\in K.$
     \end{definition }
     Obviously, that if $u_{1},...,u_{n}$ is a bases of $W$ then
     from $w=\sum\limits_{i=1}^n\lambda_{i}\cdot u_{i}$ follows
      $\delta_{W}(w)=\sum\limits_{i=1}^n\lambda_{i}^{\delta}\cdot u_{i}$.
   \begin{definition }

          Let $F_{1}$ and $F_{2}$ be Lie algebras over a field $K$,
       $\delta$ be an automorphism of $K$ and $\varphi :F_{1}\longrightarrow F_{2}$ be a ring
       homomorphism of these algebras. A pair $(\delta,\,\varphi)$
       is called semi-homomorphism from $F_{1}$ to $F_{2}$ if
       $$
           \varphi(\lambda\cdot u)=\lambda^{\delta}\cdot
           \varphi\,(u),\;\;\mbox{for}\;\alpha\in K,\;u\in F_{1}
       $$
       \end{definition }
     \begin{proposition }
     A pair $(\delta,\,\delta_{W})$ is an semi-automorphism of $W$.
     \end{proposition }
     \begin{proof}
     It is sufficient to check that
     \begin{equation}\label{for1}
      \delta_{W}\,[u,w]=[\delta_{W}\,u,\delta_{W}\,w]
     \end{equation}
      for $u,\,w\in W$. Let $V= <v_{1},v_{2},...,v_{n},....>$ be the Hall-Shirshov
     base \cite{Bah}, \cite{Bok}. Represent $u$ and $w$ in the basis $V$.
     It is well-known that a structural constants for this
     basis are integers \cite{Bok}. The formula (\ref{for1}) is a consequence of this fact.
      \end{proof}

      We will call $\delta_{W}$ the field semi-automorphism of
   $F$ defined by a field automorphism $\delta$.
   \begin{definition }
    An automorphism $T$ of $End\,W$ is called semi-inner if it
    can be presented in the form: $T= \rho\cdot\hat{\delta}_{W}$,
    where $\hat\delta_{W}(\beta)=\delta_{W}^{-1}\cdot\beta\cdot\delta_{W},\;\beta\in
    End\,W$, $\delta_{W}$ is a field semi-automorphism of $End\,W$ and $\rho$
    is an inner automorphism of $End\,W$.
    \end{definition }
    \begin{theorem} \label{m2}
      Every automorphism of semigroup $End\,W$, where $W=W(x,y)$, is semi-inner.
      \end{theorem}
   We prove the Theorem \ref{m2} in several steps.

  I. It is known \cite{Cohn} that $Aut\;W\simeq
GL_{2}\;{(K)}$. Therefore
 \begin{equation}\label{for2}
 Aut\,Aut\,W\simeq Aut\,GL_{2}\;{(K)}
 \end{equation}
 However the description of $Aut\,GL_{2}\;{(K)}$ is known \cite{Die}.
 Namely, if $S \in Aut\,GL_{2}(K)$, then there are two
possibilities:

1. $S(\beta) =\chi (\beta)\cdot\tau
\cdot\,f(\beta)\cdot\tau^{-1}$, where $\beta,\,\tau \in
Aut\;GL_{2}$, $\chi : Aut\,W \longrightarrow K^{*}$ is a
multiplicative homomorphism, $f$ is an automorphism of the field
$K$ and if $\beta = \left (\begin{array}{cc} a_{11}&
a_{12}\\a_{21}& a_{22} \end{array}\right )\in GL_{2}(K)$, then
$f(\beta) = \left (\begin{array}{cc} f(a_{11})&
f(a_{12})\\f(a_{21})& f(a_{22}) \end{array}\right )$.

2.  $S(\beta) =\chi (\beta)\cdot\tau \cdot f\,(\widetilde{\beta})
\cdot\tau^{-1}$, where $\chi$, $\tau$, $\beta$ and $f$ are as
above and $\widetilde{\beta}$ is the transformation
contragradient to $\beta$.

Note that if $\beta\in GL_{2}(K)$ then $\widetilde \beta =
\delta^{-1}\cdot\beta \cdot \delta$ for some $\delta\in GL_{2}$.
This follows from the conjugacy of Jordan forms of the matrix
$\beta$ and its  contragradient matrix $\widetilde \beta$
 in $GL_{2}(K)$. Therefore, up to conjugacy, it is sufficient to consider only the case
 1.
  By virtue of the isomorphism (\ref{for2}) we conclude that if $T\in
  Aut\,Aut\,W$ and $\xi\in Aut\,W$, then there exists $\tau\in
  GL_{2}(K)$ and $\gamma\in Aut\,K$ such that
  $$
  T(\xi)=\chi(\xi)\cdot\tau^{-1}\cdot\gamma_{W}^{-1}\cdot\xi\cdot\gamma_{W}\cdot
  \tau.
  $$
  Putting
  $$
  T_{1}(\xi)=\gamma_{W}\cdot\tau\cdot
  T(\xi)\cdot\tau^{-1}\cdot\gamma_{W}^{-1},
  $$
   we obtain
    $$
     T_{1}(\xi)= \chi(\xi)\cdot\xi ,\,\,\xi\in Aut\,W.
    $$
     \begin{definition }
   Let $T\in\,Aut\,End\,W$ and the restriction $T|_{Aut\,W}=T_{1}$.
   Automorphism $T$ is called a characteristic
   automorphism of $End\,W$ if  $T_{1}=\chi_{1} \cdot I$, where $I$ is
   the identical automorphism on $Aut\,W$ and $\chi_{1}$ is a multiplicative
   homomorphism $Aut\,W$ on $K^{*}$.
    \end{definition }

    Let $T$ be an automorphism of the semigroup of endomorphisms of
    $End\,W$. Then, obviously, $T \in Aut\,(AutF)$ and, replacing $T$
    by $T_{1}$, as above we can regard $T$, up to conjugacy, as a characteristic
    automorphism on $Aut\;W$.

    II. Now consider a $L$-endomorphism of $W$.
    \begin{definition }
 An endomorphism $\sigma \in End\,W$ is called L-endomorphism of $W$ if
 \begin{equation*}
  \begin{array}{cc}
 \sigma(x)=a_{11}x +a_{12}y \\
  \sigma(y)=a_{21}x +a_{22}y, \\
 \end{array}
 \end{equation*}
 where $a_{ij}\in K$.
  \end{definition }
   Denote by $D$ the semigroup of L-endomorphisms of $W$.  There
   exists the natural isomorphism $\varphi:\, D\longrightarrow
   Mp_{2}(K)$, where $Mp_{2}(K)$ is the semigroup of $2\times 2$ matrix
   over the field $K$.
   \begin{proposition } \label{pra1}
     Let $T \in Aut\,(End\,W)$ be a characteristic automorphism of
 $End\,W$, where $W$ is a free $2$-generated Lie algebra over a infinite field $K$,
 and $\sigma\in D$. Then there exists $\gamma\in Aut\,K$ and $\tau\in
 Aut\,W$ such that
 \begin{equation}\label{for0}
 T(\sigma)=\tau\cdot\gamma_{W}\cdot\sigma\cdot\gamma_{W}^{-1}\cdot\tau^{-1}
  \end{equation}
 \end{proposition }
 \begin{proof}
 1. We check first that our assertion holds for the idempotent
  $\sigma\in End\,W$: $ \sigma(x)=x$ and $\sigma(y)=0$.
 Let
  \begin{equation} \label{for3}
 \begin{array}{cc}
  T\,\sigma (x)= a_{11}\,x+a_{12}\,y + \sum\limits_{n=3}a_{1n}\,u_{n}(x,y)\\
   T\,\sigma (y)= a_{21}\,x+a_{22}\,y + \sum\limits_{n=3}a_{2n}\,u_{n}(x,y)\\
   \end{array}
  \end{equation}
  be the representations of $T\,\sigma (x)$ and $T\,\sigma (y)$
  relative to same basis of $W$, $a_{ij}\in K$, and $deg\,u_{n}\geq 2$ for $n\geq 3$.
  Denote by $cont_{x}\,(u)$  (or $cont_{y}\,(u) )$
  a number of appearances of $x$ (or $y$) in the word $u$.
  Now consider the automorphism $\varphi_{\lambda}$ of $W_{2}$:
 $\varphi_{\lambda}(x)=x,\;
\varphi_{\lambda}(y)=\lambda\,y,\;\lambda\neq 0,1$. Then
$\varphi_{\lambda}\cdot\sigma= \sigma$. Therefore,
$T\;\varphi_{\lambda}\cdot T\;\sigma= T\;\sigma$, i.e.
\begin{equation}\label{fort}
 k\cdot\varphi_{\lambda}\cdot T\,\sigma= T\,\sigma,
 \end{equation}
where $k=\chi(\varphi_{\lambda})\neq0$. Substituting (\ref{for3})
into (\ref{fort}) and comparing the corresponding coordinates we
obtain
 $$
k\,a_{1n}\,\lambda^{r}= a_{2n}, \;\;\mbox{where}\;
r=cont_{y}\,(u_{n}),\;n\geq 2.
 $$
  Since the field $K$ is infinite, $a_{1n}=a_{2n}=0$ for $n\geq
  2$. Hence, $T\,\sigma$ is $L$-endomorphism of $W$.

 2. Now we consider the general case.  Let
$\sigma=\left(\begin{array}{cc}
  a_{11}& a_{12}\\
  a_{21}& a_{22}\\
   \end{array}\right)$
  be an $L$-endomorphism of $W$. If $det\,\sigma\neq 0$ then $\sigma$
is automorphism of $W$ and, therefore, $T\,\sigma$ is
$L$-endomorphism.
 If $det\,\sigma=0$, then there exists $d\in K$ such that
   $$
    \begin{array}{lcc}
   \sigma(x)=a_{11}\,x +a_{12}\,y \\
    \sigma(y)=d\,a_{11}\,x +d\,a_{12}\,y ,\\
   \end{array}
$$
 The eigenvalues of matrix
  $\sigma=\left(\begin{array}{cc}
  a_{11}& a_{12}\\
  d\,a_{11}& d\,a_{12}\\
   \end{array}\right)$ are
  $\lambda_{1}=a+db,\;
 \lambda_{2}=0$. Therefore the possible Jordan forms of the matrix $\sigma$ are
 $\left(\begin{array}{cc} \lambda_{1}& 0\\0& 0
\end{array}\right )$ or $\left(\begin{array}{cc} 0& 0\\0&
\lambda_{1} \end{array}\right )$.

Consider only the first case. Let
$\delta(x)=\lambda_{1}\,x,\;\delta(y)=0,\;\lambda_{1}\neq 0$.
 If $\varphi_{\lambda_{1}}(x)=\lambda_{1}\,x,\;\varphi_{\lambda_{1}}(y)=y$
and $\sigma_{1}(x)=x,\;\sigma_{1}(y)=0$ then
$\delta=\varphi_{\lambda_{1}}\,\sigma_{1}$. Therefore,
  $$
  T\,\delta=T\,\varphi_{\lambda_{1}}\cdot T\,\delta_{1}=\chi(\varphi_{\lambda_{1}})\cdot
   \varphi_{\lambda_{1}}\cdot T\,\sigma_{1}.
  $$
  Since $\sigma_{1}$ is idempotent of $W$,
  $T\,\sigma_{1}$ is an L-endomorphism of $W$ (see point 1),
  $T\,\delta$ is also an L-endomorphism, so $T\,\sigma$ is
 $L$-endomorphism as well. We had proved that $T\;D=D$ for all $T\in Aut End \;W$.
 Now you need in the following result
  \begin{proposition }\cite{Lam}\label{pram1}
  Let $M_{n}(K)$ be a full matrix semigroup over a field $K$ and
  $f:M_{n}(K)\longrightarrow M_{n}(K)$ is semigroup automorphisms of $M_{n}(K)$.
   Then for some $\gamma\in Aut\;K$ and  for some $\Phi\in GL_{n}(K)$ holds
     \begin{equation*}
     f(A)= \Phi\cdot\gamma(A)\cdot\Phi^{-1},\,\; A\in M_{n}(K),
      \end{equation*}
      where if $A=(a_{ij})$ then $\gamma(A)=(\gamma(a_{ij}))$.
    \end{proposition }
     Now the formula (\ref{for0}) follows from Proposition \ref{pram1}.
     The proof of the theorem is complete.
  \end{proof}
    \begin{remark}\label{rem}
   From Proposition \ref{pra1} follows that any characteristic automorphism $T$ of $End\,W$
  is the identical automorphism on the semigroup $D$ of
  $L$-endomorphisms of $W$ up to conjugacy.
    \end{remark}
   \begin{proposition }\label{pra2}
  If $T$ is a characteristic automorphism of $End\,W$ and $\sigma\in
  End\,W$
  then $T\,\sigma= \tau_{\lambda}^{-1}\,\sigma\,\tau_{\lambda}$,
  where $\tau_{\lambda}(x)=\lambda\,x, \,
  \tau_{\lambda}(y)=\lambda\,y,\,\lambda\neq 0\in K $.
  \end{proposition }
  \begin{proof}
  From now on we assume that any characteristic automorphism of $End\,W$
  is the identical automorphism on the semigroup $L$-endomorphisms of
  $W$ (see. Remark \ref{rem}).

   1. Denote by $d(p)$ an endomorphism of $W$ such that $d(p)(x)=p,\;d(p)(y)=y,\;p\in W$.
 We would like to show that
  \begin{equation}\label{for4}
  T\,d(p)= \tau_{\lambda}^{-1}\,d(p)\,\tau_{\lambda}.
   \end{equation}
  (i) We assume first that $p=[x,y]$, i.e.
  $d(p)(x)=[x,y],\;d(p)(y)=y$. Let
  $$
  T\,d(p)(x)=\alpha\,[x,y]+A(x,y),
  $$
   where $A(x,y)\in W$. Let
   $\varphi_{\lambda}(x)=\lambda\cdot x,\,\varphi_{\lambda}(y)=y$.
   Then $d([x,y])\cdot\varphi_{\lambda}=\varphi_{\lambda}\cdot
   d([x,y])$. Therefore
      \begin{equation}\label{for5}
   \lambda\cdot A(x,y)=A(\lambda\cdot x,y).
     \end{equation}
   Let $A(x,y)=\sum\limits_{i}m_{i}$ be a linear combination of Lie monomials $m_{i}$ from $W$.
   From the formula (\ref{for5}) follows that $cont_{x}(m_{i})=1$ for all $i$.
    Now consider two endomorphisms $\varphi_{1},\,\varphi_{2}$
   of $W$: $\varphi_{1}(x)=\frac{1}{k}\,x,\,\varphi_{1}(y)=k\,y$
   and $\varphi_{2}(x)=x,\,\varphi_{2}(y)=k\,y$ for $k\neq 0$. Then
   $\varphi_{1}\cdot d([x,y])=d([x,y])\cdot\varphi_{2}$.
   Therefore,
   $(\varphi_{1}\cdot T\,d([x,y])(x)=(T\,d([x,y])\cdot\varphi_{2})(x)$.
   Finally we have $\varphi_{1}\,A(x,y)=A(x,y)$. Since
   $cont_{x}(m_{i})=1$ for all $i$ we obtain $A(x,y)=\beta \,[x,y]$,
   i.e., $T\,d(p)(x)=\lambda\,[x,y]$, ($\lambda=\alpha+\beta$).

   Now we consider the endomorphism $\varphi$ of $W$: $\varphi(x)=\varphi(y)=y$.
   Since $d(p)\cdot\varphi=\varphi$,
  $(T\,d(p))\,\varphi=\varphi$. Therefore $T\,d(p)(y)=y$. We have that
  the equality (\ref{for4}) is fulfilled in this case .

   (ii) Let $p$ be a right-normalized monomial of length $n$ from $W$.
   We prove the validity of (\ref{for4}) using induction on $n$. We have the basis of
   induction in the part (i). Two cases are possible: $p=[v,y]$
   or $p=[v,x]$, where $deg\,(v)<n$. We investigate only the first case.

   Consider the endomorphisms $d(v)$ and $d([x,y])$. It is
   obvious that  $d(v)\cdot d([x,y]= d(p)$. Therefore,
    \begin{equation}\label{for6}
     T\,d(v)\cdot T\,d([x,y]=T\,d(p).
    \end{equation}
  According to (i), we have
  $T\,d([x,y])=\tau_{\lambda}^{-1}\,d([x,y]\,\tau_{\lambda}$ and
  by induction $T\,d(v)=\tau_{\lambda}^{-1}\,d(v)\,\tau_{\lambda}$.
  Substituting the last equalities in (\ref{for6}) we obtain the desired.

  (iii) Now consider an endomorphism $\varphi$:
  $\varphi\,(x)=v_{1}$, $\varphi\,(y)=v_{2}$, where $v_{1},\;v_{2}$ are
  right-normalized monomials.
   We have
       \begin{equation} \label{for7}
       T\,d(\varphi\,(x))=T\,(\varphi\,d(x))=(T\,\varphi)d(x)=d(T\,\varphi(x))
    \end{equation}
    However from (ii) follows that
    \begin{equation}\label{for8}
    T\,d(\varphi(x))=
    \tau_{\lambda}^{-1}\,d(\varphi(x))\,\tau_{\lambda}
    =d({\tau_{\lambda}^{-1}\,\varphi\,\tau_{\lambda}}(x)).
    \end{equation}
   Comparing (\ref{for7}) and (\ref{for8}) we obtain
    $T\,\varphi(x)={\tau_{\lambda}^{-1}\,\varphi\,\tau_{\lambda}}(x)$.
    In a similar way we may obtain
    $T\,\varphi(y)={\tau_{\lambda}^{-1}\,\varphi(y)\,\tau_{\lambda}}$,
    i.e., $T\,\varphi={\tau_{\lambda}^{-1}\,\varphi\,\tau_{\lambda}}$.

    (iv) Let $p=v_{1}+v_{2}$, where $v_{1},\,v_{2}$ are right-normalized
    monomials. Consider the endomorphisms $\delta$ and $\varphi$:
    $\delta\,(x)=x+y,\;\delta\,(y)=y $ and
    $\varphi\,(x)=v_{1},\;\varphi\,(y)=v_{2}$. Then $d(p)=
    \varphi\cdot\delta$. From (iii) follows that
    $T\,\varphi={\tau_{\lambda}^{-1}\cdot\varphi\cdot\tau_{\lambda}}$
     and $T\,\delta =\delta$. This proves our assertion for two
     summands. Easy induction on number of summands of $p$ leads
     to the formula (\ref{for4}) in the general case.

     2. Let $\sigma \in End\,W$ such that $\sigma(x)=p_{1},\;
     \sigma(y)=p_{2}$, where $p_{1},\;p_{2}\in W.$. The same arguments
    as in  part (iii) prove the Proposition \ref{pra2}.
      \end{proof}

      {\bf Proof of Theorem \ref{m2}}. Let $T$ be an automorphism of $End\,W$ and $T_{1}$ be the
       restriction of $T$ on $Aut\,W$, i.e. $T_{1}=T\mid_{Aut\,W}$. We have for all
        $\beta\in Aut\;W$:
        $T_{1}(\beta)=\tau^{-1}\cdot\delta_{W}^{-1}
        \cdot\beta\cdot\delta_{W}\cdot\tau,\;\beta\in End\,W$, where $\delta_{W}$
        is a field semi-automorphism of $End\,W$, $\tau\in Aut\,W$.

        Let $T_{2}$ be mapping of $End\;W$ the following form:
         $T_{2}(\beta)=\tau^{-1}\cdot\delta_{W}^{-1}
        \cdot\beta\cdot\delta_{W}\cdot\tau,\;\beta\in End\,W$, for
        the same $\delta_{W}$, $\tau$ as above and for all $\beta\in End\;W$.
        It is easy to check that $T_{2}$ is an automorphism of
        $End\,W$.

       It is clear that $T_{2}\mid_{Aut\,W}=T_{1}=T\mid_{Aut\,W}$.
       Therefore,
       $T\cdot T_{2}^{-1}\,(\beta)= \beta,\,\beta\in Aut\,W$. By Proposition \ref{pra2} we have
      $T\cdot T_{2}^{-1}=\hat{\tau_{\lambda}}$. Hence, $T=\hat{\tau_{\lambda}}\cdot T_{2}
      =\hat{\tau_{\lambda}}\cdot\hat{\tau}\cdot\hat{\delta_{W}}=\phi\cdot\hat{\delta_{W}},$
      where $\phi=\hat{\tau_{\lambda}}\cdot\hat{\tau}$ is the inner automorphism
      of $End\,W$ and $\hat{\delta_{W}}$ is as above. The proof of Theorem \ref{m2} is complete.
       \subsection{Semi-inner automorphisms of the category of free Lie algebras}
       Consider the category $\bf \Theta=Lie-K $ of free Lie
       algebras over a field $K$ with the semimorphisms defined in
       \ref{dl1}.
       We can to define semi-inner automorphisms of the category $\bf \Theta^{0}$
        on the same way as in the case of the category of free modules.
      \begin{definition }
    An automorphism $\varphi$ of the category $\bf\Theta^{0}$ is
  called semi-inner if there exists  $\sigma\in Aut\,K$ and a
  collection $(\sigma,s)=\{(\sigma,s_{X})|s_{X}:F\longrightarrow F,\;F\in Ob\,\bf\Theta \}$
  of semi-isomorphisms such that $(\sigma, s) :1_{\bf \Theta^{0}}\longrightarrow \varphi$ is
  semi-isomorphism of functors, i.e., for every object $F\in Ob\,\bf \Theta$
 there exists semi-isomorphism $(\sigma, s_{X}): F\longrightarrow
 F$  and for each $\nu: F_{1}\longrightarrow F_{2}$ the
 following diagram is commutative:

    $$\CD F_{1} @>s_{F_{1}} >> \varphi(F_{1})\\ @V\nu VV @VV\varphi (\nu) V\\
F_{2} @>s_{F_{2}} >> \varphi(F_{2})\endCD$$

   \end{definition }
  Denote by $\hat{\sigma}$  an automorphism of the category of the
 free $K$-modules of the following form: \par\medskip
 1. \,$\hat{\sigma}(KX)=KX, \,KX\in Ob\;\bf C^{0}$, i.e., $\hat{\sigma}$ does not change objects
  of the category of the free $K$-modules. \par\medskip
 2. \,$\hat{\sigma}(\nu)=\sigma_{Y}\nu\,\sigma_{X}^{-1}:KX\longrightarrow KY $
  for all $\nu:KX\longrightarrow KY;\; KX,KY\in Ob\;\bf C^{0}$

  It is clear that $\hat{\sigma}$ is
  a semi-inner automorphism of the category $\bf \Theta^{0}$.

      Now, we shall describe the group of automorphisms of the category of
   free Lie algebras over an infinite field $K$.
 \begin{theorem}\cite{Mash1}\label{m3}
  Let $K$ be a infinite field. Then every automorphism of the category
  $\bf \Theta^{0}=(Lie-K)^{0}$ is semi-inner.
 \end{theorem}
 \begin{proof}
A free Lie algebra over field is hopfian \cite{Bok}.
 It is clear that if $\varphi\in Aut\;\bf \Theta^{0}$ then
  for any one-dimensional Lie algebra $Kx_{0}$ we have $\varphi(Kx_{0})\simeq Kx_{0}$.
 By Theorem \ref{t1} every automorphism of the category $\bf
 \Theta^{0}$ is special. By Proposition \ref {pr0} an automorphism
 $\varphi$ can be represented in the form
 $\varphi=\varphi_{1}\cdot\varphi_{2}$, where $\varphi_{1}$ is an inner
 automorphism of $\bf \Theta^{0}$ and $\varphi_{2}$ is an automorphism which does not change
  the objects of this category. Now, without loss of
  generality, we can assume that $\varphi$ itself does not change the objects of
  the category $\bf \Theta^{0}$. From this follows that $\varphi$
  induces  an automorphism $\varphi _{W}$ of $End\;W$
  for all $W\in \bf \Theta$. Take 2-generated Lie algebra
  $W^{0}=W(x,y)$. By Theorem \ref {m2} the automorphism   $\varphi _{W^{0}}$
  is semi-inner.  This automorphism is defined by $\delta\in Aut\,K$
  and $\tau_{W}\in Aut\,W^{0}$.

   Making use of the pair $(\delta,\tau_{W})$ we build an inner automorphism
  of the category $\bf \Theta^{0}$. Let $W\in \Theta$. For any $W\in \Theta$
  set $(\delta,\tau)_{W}=(\delta,\tau_{W})$ if $W=W^{0}$ and
   $(\delta,\delta_{W})$ otherwise. The collection
   $(\delta,\tau)_{W}$ of automorphisms of this category defines an inner automorphism
   $\phi$ of $\Theta^{0}$. Denote by $\varphi_{1}=\phi^{-1}\varphi$. Then $\varphi_{1}$
   acts identical on $End\,W^{0}$.

  Let $\nu_{0}:W^{0}\longrightarrow Kx_{0}=W_{0}$ be a homomorphism such that
    $\nu_{0}(x)=\nu_{0}(y)=x_{0}$ and $\rho: W^{0}\longrightarrow W^{0}$
    be a homomorphism satisfying $\rho(x)=y,\,\rho(y)=x$.
    From $\nu_{0}\rho=\nu_{0}$ follows
    $\varphi(\nu_{0})\rho=\varphi(\nu_{0})$. Therefore
    $\varphi_{1}(\nu_{0}(x))=\varphi_{1}(\nu_{0}(y))=\alpha x$ for some $\alpha\neq 0$
    from $K$.

 Consider the collection of automorphisms $f_W:W\longrightarrow W$, $W=W(X)\in
 \Theta$ defined in the following way:
  if $W=W^{0}$, then $f_W^{0}(x)=\alpha x_{0}$ and if $W\neq W^{0}$, then
  $f_{W}(x)=x$ for all $x\in X$. This collection of automorphisms defines an inner automorphism
  $\hat{f}$ of the category $\bf \Theta^{0}$ which does not change objects of
   $\bf \Theta^{0}$ and
   $\hat{f}(\nu)=f_{Kx_{0}}\nu f_{KX}^{-1}$ for every morphism $\nu\in Mor\; \bf \Theta^{0}$. It is clear that
  $\hat{f}(\nu_{0})=\varphi_{1}(\nu_{0})$.
  By Theorem \ref{t2} the automorphism $\psi=\hat{f}^{-1}\varphi_{1}$ is an inner
  automorphism of the category $\bf \Theta^{0}$ and thus, $\varphi=\phi\hat{f}^{-1}\psi$
  is semi-inner. The proof is complete.
  \end{proof}


\begin{thebibliography} {99}
   \bibitem{Bah} Yu. Bahturin, \textit{Identical relations in Lie
   algebras}, VNU Science Press, Utrecht, {1987}.
   \bibitem{Ber} A. Berzins, B. Plotkin, E. Plotkin,
   \textit{Algebric geometry in varieties of algebras with the
   given algebra of constants}, Journal of Math. Sciences, \textbf
   {102:3}, (2000), 4039-4070.

   \bibitem{Bok} L. Bokut, G. Kukin, \textit{Algorithmic and
   combinatorial algebra}, Kluver, (1994).

   \bibitem{Cohn} P. Cohn, \textit{Subalgebras of free associative algebras}, Proc. London
   Math. Soc, \textbf {14}, (1985), 618-632.

   \bibitem{Die}  J. Dieudonne, \textit{On the automorphisms of the
   classical groups}, Memoirs Amer. Math. Soc., \textbf{2}, (1951), 1-95.

    \bibitem{For} E. Formanek, \textit{A question of B. Plotkin about
    the semigroup of endomorpjsms of a free group}, Proc. American Math.
    Soc., \textbf{130}, (2001), 935-937.

    \bibitem{Lam} M. Jodiet, T. Lam, \textit{Multiplicative maps of matrix
    semigroups},
     Arch. Math., {\textbf 20}, (1969), 10-16.

    \bibitem{Mash0} G. Mashevitzky, B. Schein, \textit{Automorphisms
    of the endomorphism semigroup of a free monoid or a free
    semigroup}, to appear.

     \bibitem{Mash1} G. Mashevitzky, B. Plotkin, E. Plotkin, \textit{Automorphisms
      of the category of free Lie algebras}, to appear.

     \bibitem{Mash2} G. Mashevitzky, B. Plotkin, E. Plotkin, \textit{Automorphisms
      of the category of free algebras of varieties},  Electron. Res. Announs.
      Amer. Math. Soc., {\textbf 8}, (2002), 1-10.

    \bibitem{Plot}  B. Plotkin, \textit{Seven lectures in universal algebraic geometry},
      Preprint, Hebrew University, Jerusalem, (2000).

    \bibitem{Shir} A. Shirshov, \textit{Subalgebras of free Lie
    algebras}, Uspekhi Mat. Nauk, \textbf{8}, (1953), 173-176.

    \bibitem{Hal} M. Hall,  \textit{A basis for free Lie ring and higher
    commutators in free groups}, Proc. London Math. Soc., \textbf{1}, (1950), 575-581.

     \bibitem{Zhit} G. Zhitomirskii, \textit{Autoequevalences of
     categories of free algebras of varieties}, to appear.

 \end{thebibliography}
\end{document}